\documentclass[]{article}
\usepackage{geometry}
\usepackage{moreverb}

\usepackage{amssymb,amsmath,amsthm,amsfonts}
\usepackage{enumerate}
\usepackage{xspace}
\usepackage{tikz-cd}
\usepackage{graphicx}

\newcommand{\ms}{\eqref{eq:ms}\xspace}

\newcommand{\T}{\mathrm{T}}

\newcommand{\rk}{\mathrm{rank \ }}

\newcommand{\fp}[2]{\frac{\partial #1}{\partial #2}}
\newcommand{\spn}[1]{\mathrm{span}\left\{#1\right\}}
\newcommand{\ti}[1]{\tilde{#1}}

\newtheorem{definition}{Definition}
\newtheorem{remark}{Remark}
\newtheorem{theorem}{Theorem}
\newtheorem{corollary}{Corollary}
\newtheorem{proposition}{Proposition}

\begin{document}

\title{Input-output linearization and decoupling of mechanical control systems}

\author{Marcin Nowicki,\\ Institute of Automatic Control and Robotics, Poznan University of Technology,\\ Pozna\'n, Poland, marcin.nowicki@put.poznan.pl \vspace{0.5cm}\\
Witold Respondek,\\ Institute of Automatic Control, Lodz University of Technology, \L{}\'od\'z, Poland\\
Emeritus Professor at Normandie Universit\'e, INSA de Rouen, Saint-Etienne-du-Rouvray, France\\ witold.respondek@insa-rouen.fr
}

%
%
%
%
%
%
%




\maketitle


\section*{Abstract}
{In this work, we present the problem of simultaneous input-output feedback linearization and decoupling (non-interacting) for mechanical control systems with outputs. We show that the natural requirement of preserving mechanical structure of the system and of transformations imposes supplementary conditions when compared to the classical solution of the same problem for general control systems. These conditions can be expressed using objects on the configuration space only. We illustrate our results with several examples of mechanical control systems. }

\section{Introduction}

In this paper, we formulate and solve the problem of simultaneous input-output feedback linearization and decoupling of mechanical control systems, shortly MIOLD-problem. Compared to our previous results \cite{aut,tac,lh,NR}, here we consider mechanical systems with outputs, denoted \ms, and give conditions for the output functions to form (partially or fully) linearizing outputs. We thus relate our problem to the classical works \cite{ida,nij,Isidori,cirt} on the problem of input-output decoupling (i.e. non-interacting control) and input-output linearization, which is formalized as follows. Consider a square nonlinear control system with outputs of the form
	\begin{align}
		\label{eq:cs}
		\dot{z}= F(z)+\sum_{r=1}^{m}G_r(z)u_r=F(z)+G(z)u, \qquad y=h(z)=(h_1(z),\ldots,h_m(z))^T,
	\end{align}
	where $z\in Z$, a manifold of dimension $N$, $u\in \mathbb{R}^m$, $
	y \in \mathbb{R}^m$. Find a static invertible feedback $	u_r= \alpha^r(z)+\sum_{s=1}^{m}\beta^r_s(z)\tilde{u}_s$ such that for the closed-loop system the $i$-th input $\ti u_i$  affects only the corresponding $i$-th output $y_i=h_i(z)$ and does not affect the other outputs $y_j$, for $j\neq i$.
	{A well known result, see e.g. \cite{ida,nij,Isidori} asserts that the problem of input-output linearization and decoupling is solvable if and only if $\rk{\mathbb{D}(z)}=m$, where $\mathbb{D}=(L_{G_r}L_F^{\rho_\ell-1}h_\ell)$, $1\leq r,\ell\leq m$, is the decoupling matrix and $\rho_\ell$ are relative degrees, see Section \ref{sec:PF} for all definitions and the corresponding normal  form.}
	
{In this paper, we will consider the input-output linearization and decoupling problem in the case, where system \eqref{eq:cs} is mechanical, that is, the state $z=(x,v)$ consists of configurations and velocities, and both, the drift $F$ and the control vector fields $G_r$ exhibit a mechanical structure, see Section \ref{sec:PF} for a precise definitions and Section \ref{sec:MS} for additional information on mechanical control systems.}
We shall consider the following questions:
\begin{enumerate}[Q1.]
	\item For a mechanical control system with outputs, how to render its input-output relation both linear and input-output decoupled in a way that preserves the mechanical structure of the system?
	\item How does the problem of simultaneous input-output linearization and decoupling for mechanical systems differ from that problem for general control systems?
	\item Is it possible to formulate conditions using objects on the configuration manifold only (instead of those on the whole state-space)?
	\item What are properties of both, the observable and the unobserved dynamics associated with the above MIOLD-problem?
\end{enumerate}

In Section \ref{sec:PF}, {we describe the class of mechanical control systems and} formulate the problem of mechanical input-output linearization and decoupling, also known as non-interacting control problem. In Section \ref{sec:MS}, we present foundations of geometric view of mechanical control systems and make connections between the class of mechanical control systems and conservative Lagrangian control systems. In Section \ref{sec:IO}, we give the main result of the paper, which is a geometric solution of the MIOLD-problem, i.e. we give conditions under which {for} a mechanical system {we can render, via static invertible feedback,}  the input-output relation both linear and decoupled (non-interactive), and simultaneously mechanical. Then, we state conditions for mechanical feedback linearization for systems without outputs (for which we add virtual outputs) via feedback transformations that preserve the mechanical structure of the system, thus relating the MIOLD-problem and the results of this paper to those of references \cite{aut,tac}  .

\section{Problem formulation}
\label{sec:PF}
{We introduce a minimal amount of objects and terminology necessary to formulate the result. For more on geometric approach to mechanical systems see \cite{bullo,bloch}, geometric control theory \cite{nij,Isidori}, and identifying mechanical among general control systems \cite{rr}.} 
Consider a mechanical control system together with output, that depends on the configurations only, of the form
\begin{align}
	\begin{split}
		\dot{x}&=v\\
		\dot{v}&=-v^T\Gamma(x)v+e(x)+\sum_{r=1}^{m}g_r(x) u_r,\\
		y&=h(x) 
	\end{split}
	\label{eq:ms}
	\tag{$\mathcal{MSO}$}
\end{align}
where  $x =\left(x^1,\ldots,x^n \right)^T $ are local coordinates on the configuration manifold $Q$ of the system, and $v =\left(v^1,\ldots,v^n \right)^T $ are velocities and thus the pair $(x,v)$ forms coordinates on the tangent bundle $\T Q$, and $\Gamma(x)$ is the matrix of Christoffel symbols of the second kind $\Gamma^i_{jk}(x)$, the vector fields $e(x)=( e^1(x),\ldots,e^n(x))^T$ and $g_r(x)=(g_r^1(x),\ldots,g_r^n(x))^T$ correspond to, respectively, uncontrolled and controlled action on the system. The output $y=(y_1,\ldots,y_m)^T$ depends on configurations $x\in Q$ only, that is, $y=h(x)=(h_1(x),\ldots,h_m(x))^T$, where $h:Q\rightarrow\mathbb{R}^m$ is a smooth mapping. The {considered} system is square (equal numbers of scalar inputs and scalar outputs). A generalization to the case of more inputs than outputs is straightforward, while in the case of less inputs than outputs the non-interacting property cannot be achieved.

Throughout all objects are assumed to be smooth and the word smooth means $C^{\infty}$-smooth.	The tensor summation convention is assumed throughout, i.e. any expression containing a repeated index (upper and lower) implies summation over that index (up to $n$), e.g. 	$\omega_iX^i=\sum_{i=1}^n\omega_iX^i$; we will not, however, apply that convention to controls $u_r$ and control vector fields $g_r$, where the summation goes up to $m$.

Note that, the mechanical control system \ms evolves on the tangent bundle $\T Q$, consequently  $e^i(x)\fp{}{v^i}$ and $g_r^i(x)\fp{}{v^i}$ are vector fields on $\T Q$. Actually, those are \textit{vertical lifts} of vector fields on $Q$, i.e. $e(x)=e^i(x)\fp{}{x^i}$ and $g_r(x)=g_r^i(x)\fp{}{x^i}$, respectively, and define on $Q$ the \textit{virtual system} $\dot{x}=e(x)+\sum_{r=1}^{m}g_r(x)u_r$. We will use those vector fields on $Q$ repeatedly, in particular to formulate conditions using objects defined on the manifold $Q$ only. See e.g. \cite{bullo} for more on vertical lift of vector fields and a coordinate-free description of mechanical control systems, here we intend to work in coordinates.

Consider the group $G_{MF}$ of mechanical feedback transformations of systems \ms given by:

\begin{enumerate}[(i)] 
	\item mechanical changes of coordinates $\Phi:\  \T Q \rightarrow  \T \tilde{Q}$ of the form:
	\begin{align}
		\label{eq:diff}
		\begin{split}
			\left(x,v \right)  &\mapsto \left( \ti{x},\ti{v}\right) =\Phi(x,v)=\left( \phi(x),\fp{\phi}{x}(x)v \right), 
		\end{split}
	\end{align}
	called mechanical diffeomorphisms,	where $\phi:Q\rightarrow \ti{Q}$ is a diffeomorphism and {$\fp{\phi}{x}=\left(\fp{\phi^i}{x^j}\right)$} its Jacobi matrix;
	
	\item mechanical static feedback transformations, denoted $\left(\alpha,\beta,\gamma \right) $, of the form
	\begin{align}
		\label{eq:feed}
		u_r=v^T\gamma^r(x) v+\alpha^r(x)+\sum_{s=1}^{m}\beta^r_s(x)\tilde{u}_s,
	\end{align}
	where $\gamma^r_{jk},\alpha^r,\beta^r$ are smooth functions on $Q$ satisfying $\gamma^r_{jk}=\gamma^r_{kj}$ (thus transforming the Christoffel symbols), and the matrix $\beta(\cdot)$ is invertible.
\end{enumerate}
Notice that elements of the group $G_{MF}$ of mechanical feedback transformations preserve the mechanical structure of systems. Indeed, a mechanical diffeomorphism $\Phi(x,v)=\left( \phi(x),\fp{\phi}{x}(x)v \right)$, together with a mechanical feedback transformation $(\alpha,\beta,\gamma)$, map the system \ms into the mechanical system
\begin{align}
	\begin{split}
		\dot{\ti x}&=\ti v\\
		\dot{\ti v}&=-\ti v^T\ti \Gamma(\ti x)\ti v+\ti e(\ti x)+\sum_{s=1}^{m}\ti g_s(\ti x) \ti u_s,\\
		y&=\ti h(\ti x) 
	\end{split}
	\label{eq:mst}
	\tag{$\widetilde{\mathcal{MSO}}$}
\end{align}
where $\tilde{e}=\phi_*e$, $\ti g_s=\phi_*(\beta^r_sg_r)$, $\tilde{h}=\ti \phi^*h$, with $\ti \phi=\phi^{-1}$, and $\ti \Gamma(\ti x)$ are the Christoffel symbols $\Gamma^i_{jk}-g^i_r\gamma^r_{jk}$ expressed in $\ti x$-coordinates. Recall that for any vector field $f$ on $Q$, any function $h$ on $Q$, and any diffeomorphism $\phi:Q\rightarrow\ti Q$, one defines $(\phi_*f)(\ti x)=\fp{\phi}{x}(\ti \phi\left(\ti x)\right)f(\ti \phi(\ti x))$ and $\ti \phi^*h(\ti x)=h(\ti \phi(\ti x))$.

Consider the {general} square nonlinear control system with outputs \eqref{eq:cs}. Its relative degree $\rho_\ell$ of $h_\ell$ at $z_0\in Z$ is the smallest integer such that we have $L_{G_r}L_F^kh_\ell\equiv0$ in a neighborhood of $z_0$, for all $1\leq r\leq m$ and $0\leq k \leq \rho_\ell-2$ and in any neighborhood of $z_0$ there exist $z$ and $r$ such that $L_{G_r}L_F^{\rho_\ell-1}h_\ell(z)\neq 0$.
Define the decoupling matrix 
{
\begin{align*}
\mathbb{D}=(L_{G_r}L_F^{\rho_\ell-1}h_\ell), \qquad1\leq r,\ell\leq m,
\end{align*}}
and denote by $Z^0$ an open and dense subset of $Z$ such that  $\rk \mathbb{D}(z)$ is locally constant on $Z^0$. It is well known, compare \cite{nij,Isidori}, that the problems of input-output decoupling (non-interacting control) and that of input-output decoupling with simultaneous input-output linearization coincide on $Z^0$ and are solvable around $z_0\in Z^0$ if and only if $\rk \mathbb{D}(z_0)=m$. Moreover, after applying a suitable feedback $u=\alpha(z)+\beta(z)\ti u$, the system takes the form 
\begin{align}
	\label{eq:gnf}
	\begin{split}
		y_\ell=\tilde{z}^1_\ell \qquad \quad &	\dot{ \tilde{z}}_\ell^j=\ti z_\ell^{j+1} \qquad 1\leq j \leq \rho_\ell-1\\
		&	\dot{\tilde{z}}_\ell^{\rho_\ell}=\ti u_\ell\\
		&	\dot {\tilde{z}}_{m+1}=F_{m+1}(\ti z)+G_{m+1}(\ti z)\ti u,
	\end{split}
\end{align}
where $\ti z^j_\ell=L_F^{j-1}h_\ell$, with $1\leq \ell \leq m$, $1\leq j\leq \rho_\ell - 1$, and $\ti z_{m+1}$ is a vector of $N-\sum_{\ell=1}^{m}\rho_\ell$ functions that complete the $\ti z^j_\ell$'s to a local coordinate system.
If for the mechanical system \ms, expressed as $\dot z=F(z)+G(z)u$, $y=h(z)$, the condition $\rk \mathbb{D}(z)=m$ is satisfied, then the above result applies, but when  input-output decoupling the system  we may lose its mechanical structure; more precisely, the diffeomorphism $\ti z=\Psi(z)$  and feedback $u=\alpha(z)+\beta(z)\ti u$, that render the system input-output decoupled, need not be elements of the group $G_{MF}$ of mechanical feedback transformations. Our aim is to input-output decouple the system \ms and, simultaneously,  to preserve its mechanical structure, which we formalize as follows.
We say that the mechanical input-output linearization and decoupling problem, shortly MIOLD-problem, is solvable if there exist a diffeomorphism $\Phi$ of the form \eqref{eq:diff} and a feedback of the form \eqref{eq:feed} that bring the system \ms into
\begin{align}
	\tag{$\mathcal{MIOLD}$}
	\label{eq:md}
	\begin{split}
		y_\ell=\ti x_\ell^1\qquad	&\dot{\ti x}^j_\ell=\ti v_\ell^j\qquad \quad  1\leq j\leq \nu_\ell\\
		&\dot{\ti v}_\ell^j=\ti x_\ell^{j+1} \qquad \ 1\leq j\leq \nu_\ell-1\\
		&\dot{	\ti  v}^{\nu_\ell}_\ell=\ti u_\ell,\\
		&\dot{	\ti  x}_{m+1}=\ti v_{m+1}\\
		&\dot{	\ti  v}_{m+1}= -\ti v^T \ti \Gamma(\ti x)\ti v + \ti e(\ti x) +  \ti G(\ti x) \ti u,
	\end{split}
\end{align}
where $1\leq \ell\leq m$ and $\ti x_{m+1}$  is a vector that completes the $\ti x^j_\ell$'s to a local coordinate system {on $Q$}.

The normal form \eqref{eq:md} is mechanical in several ways. First, the mechanical feedback transformations \eqref{eq:diff} and \eqref{eq:feed} preserve the mechanical structure of \ms, i.e. \eqref{eq:md} is still a mechanical control system and $\ti x=\phi(x)$ maps the configurations $x$  of \ms into the configurations $\ti x$ of \eqref{eq:md} while $\fp{\phi}{x}$ maps the velocities $v$ into velocities $\ti v$. 
Second, the observable subsystem evolving on $(\ti x_o,\ti v_o)$, where $\ti x_o=(\ti x^j_\ell)$ and $\ti v_o=(\ti v^j_\ell)$, for $1\leq \ell \leq m$ and $1\leq j \leq \nu_{\ell}$,  is mechanical, linear and decoupled implying that the input-output relation of \eqref{eq:md} is linear and decoupled. More precisely, the $(\ti x_o,\ti v_o)$-subsystem consists of $m$ even-dimensional chains of integrators; this is the canonical form of linear mechanical controllable systems \cite{lms} and the outputs $y=(y_1,\ldots,y_m)^T$ are compatible with that form, i.e. they form the first state-variables (configurations) at each chain. The lengths $2\nu_\ell$ of the chains are defined by the relative degrees $\nu_{\ell}$, for $1\leq \ell\leq m$, of the virtual system $\dot{x}=e(x)+\sum_{r=1}^{m}g_r(x)u_r$, $y=h(x)$, which form the relative half-degree of \ms, see Definition \ref{def:reldeg} below for a precise formulation. Finally, the  equations of the unobserved part given by $(\ti x_{m+1},\ti v_{m+1})$, possess the mechanical structure as well, namely, they involve $\ti e(\ti x),\ti \Gamma(\ti x),\ti G(\ti x)$  that are vectors and matrices of appropriate sizes that depend nonlinearly on both $\ti x_o$ and $\ti x_{m+1}$ and, moreover, the Christoffel symbols $\ti \Gamma(\ti x)$ are multiplied by components of all velocities $\ti v = (\ti v_{o},\ti v_{m+1})$.

The following example illustrates the difference between the notion of input-output decoupling in the mechanical versus the general context.

\textbf{Example 1:}
	Consider the mechanical system with two inputs, two outputs, and 3 degrees of freedom
	\begin{align}
		\label{eq:ex11}
		\begin{split}
			y_1&=x^1\\
			\dot{x}^1&=v^1\\
			\dot{v}^1&=u_1
		\end{split}
		\begin{split}
			y_2&=x^2\\
			\dot{x}^2&=v^2\\
			\dot{v}^2&=-\Gamma^2_{22}(v^2)^2-\Gamma^2_{33}(v^3)^2+x^3+g^2_2 u_2
		\end{split}\qquad 
		\begin{split}
			\\
			\dot{x}^3&=v^3\\
			\dot{v}^3&=u_2.
		\end{split}
	\end{align}
	If $g_2^2(x)\neq 0$, then $\rho_1=\rho_2=2$, the decoupling matrix $\mathbb{D}=(L_{G_r}L_Fh_\ell)$ satisfies $\rk \mathbb{D} (z)=2$, {where $z=(x,v)$,} and the system can be transformed, via the mechanical feedback $\ti u_1=u_1$, $\ti u_2= -\Gamma^2_{22}(v^2)^2-\Gamma^2_{33}(v^3)^2+x^3+g^2_2 u_2$ into the form \eqref{eq:md}. Now, assume that $g_2^2\equiv0$ and that $\Gamma^2_{33}\neq 0$, then around any point $(x_0,v_0)$ such that $\Gamma^2_{33}(x_0)v_0^3\neq 0$ the vector relative degree is $(\rho_1,\rho_2)=(2,3)$, the decoupling matrix $\mathbb{D}$ is invertible and the system is input-output decoupable but neither the normal form \eqref{eq:gnf} nor linearizing transformations are mechanical. Indeed, the observable subsystem is of dimension $5$, the unobserved subsystem is of dimension~$1$ (so none of them is mechanical) and the variable $L_F^2h_2= -\Gamma^2_{22}(v^2)^2-\Gamma^2_{33}(v^3)^2+x^3$ mixes up positions and velocities. Finally, assume that $g_2^2\equiv0$, $\Gamma^2_{33}\equiv0$ and $\fp{\Gamma^2_{22}}{x^3}\equiv0$. Then the vector relative degree is well defined and $(\rho_1,\rho_2)=(2,4)$, the decoupling matrix $\mathbb{D}$ is invertible and the system is input-output decoupable. Notice that $\rho_1=2\nu_1=2$ and $\rho_2=2\nu_2=4$, where $(\nu_1,\nu_2)$ are relative half-degrees, see Definition \ref{def:reldeg} in Section \ref{sec:IO}. Nevertheless, the system is not mechanically input-output decoupable (unless $\Gamma^2_{22}\equiv 0$) because $L_F^2h_2= -\Gamma^2_{22}(v^2)^2+x^3$    mixes up positions and velocities (like in the previous case). If {$g^2_2\equiv0$, $\Gamma^2_{33}\equiv0$, and } $\Gamma^2_{22}\equiv0$, then the system is already linear and input-output decoupled.

\section{Mechanical control systems}
\label{sec:MS}

A class of mechanical control systems that we study is motivated by the class of conservative Lagrangian control systems (called also \textit{simple mechanical control system} in \cite{bullo}) but is more general  because  $e(x)$ can be any vector field (given by a, possibly, non potential force) and because of allowing for any symmetric (not necessarily metric) connection defining $\Gamma^i_{jk}$. 
In this section, we recall	various representations of mechanical control systems starting with Lagrangian control systems, then we give a geometric definition of mechanical control systems studied in the paper.

Consider a mechanical control system with $n$ degrees of freedom (DOF) and $m$ controls. Its Lagrangian, defined as the difference between the kinetic energy $T(x,\dot{x})$ and the potential energy $V(x)$, reads $\mathcal{L}=T(x,\dot{x})-V(x)=\frac{1}{2}\dot{x}^TM(x)\dot{x}-V(x)$, where the symmetric positive definite matrix $M(x)$ is the inertia matrix (a metric tensor) of the system. We assume that there is no energy dissipation (e.g. friction, damping) in the system and that it is subject to two kinds of external forces that are positional: external control forces $\sum\limits_{r=1}^{m}\tau_r(x)u_r$ and an uncontrolled external (not necessarily potential) force $\tau_0(x)$. The corresponding controlled Euler-Lagrange equations are
\begin{align}
	\frac{d}{dt}\fp{\mathcal{L}}{\dot{x}}-\fp{\mathcal{L}}{x}=\tau_0(x)+\sum_{r=1}^{m}\tau_r(x)u_r,\nonumber
\end{align} 
giving, with the notation $v=\dot{x}$, 
\begin{align}
	M(x)\dot{v}+C(x,v)v-P(x)=\sum_{r=1}^{m}\tau_r(x)u_r,\nonumber
\end{align}
where $C(x,v)$ is the Coriolis matrix, $P(x)=\fp{\mathcal{L}}{x}+\tau_0$ is an uncontrolled force (possibly non potential because of $\tau_0$) and $\tau_r(x)$ are external forces controlled by the controls $u_r$. Inverting the inertia matrix $M(x)$ gives $\dot{v}=-M^{-1}(x)C(x,v)v+M^{-1}(x)P(x)+M^{-1}(x)\sum\limits_{r=1}^{m}\tau_r(x)u_r$ which, in (local) coordinates $(x,v)\in\T Q$, takes the form \ms (assuming the outputs $y=h(x)$ to depend on the configurations only). The components of $M^{-1}(x)C(x,v)v$ are $\Gamma^i_{jk}(x)v^jv^k$, with $\Gamma^i_{jk}$ being the Christoffel symbols of the Levi-Civita connection corresponding to the metric $M(x)$, $e(x)=M^{-1}(x)P(x)$, and $g_r(x)=M^{-1}(x)\tau_r(x)$.

A mechanical control system on $\T Q$, with outputs, of the form \eqref{eq:ms} can be represented by a $5$-tuple $(Q,\nabla,\mathfrak{g},e,h)$ consisting of the configuration manifold $Q$ and four geometric objects defined as follows: a symmetric affine connection $\nabla$ on $Q$, an $m$-tuple of control vector fields $\mathfrak{g}=(g_1,\ldots,g_m)$ on $Q$, an uncontrolled vector field $e$ on $Q$, and a map $h:Q\rightarrow\mathbb{R}^m$.
For more on different classes of mechanical control systems see \cite{bullo,bloch, rr,rr3,rr2}.

A curve $x(t):I\rightarrow Q$, $I\subset \mathbb{R}$, is a trajectory of \eqref{eq:ms} if it satisfies the following equation
\begin{align}
	\nabla_{\dot{x}(t)}\dot{x}(t)=e\left(x(t) \right) +\sum_{r=1}^{m}g_r\left( x(t)\right) u_r,
	\label{ms0}
\end{align}
{together with the output response $y(t)=h(x(t))$. Equation \eqref{ms0} can be viewed as an equation that balances accelerations of the system, where the left hand side represents geometric accelerations (i.e. accelerations caused by the geometry of the system) and the right hand side represents accelerations caused by external actions on the system (controlled or not). System \eqref{ms0} in local coordinates $(x,v)$ on $\T Q$ takes the form of a first-order system of differential equation \eqref{eq:ms}, as given at the beginning of Section \ref{sec:PF}, to which we add the output $y=h(x)$, thus establishes a one-to-one correspondence between the differential equations of \eqref{eq:ms} and the 5-tuples.

\section{Main result}
\label{sec:IO}
The output $y=h(x)$ of \ms is given by a map $h: Q\rightarrow\mathbb{R}^m$, however the system \ms evolves on $\T Q$ so we will also interpreted $h$ (actually its pullback $\pi^*h$, where $\pi: \T Q\rightarrow Q$ is the canonical projection) as $h:\T Q\rightarrow \mathbb{R}^m$. We will denote both maps by $h$ and we use both of them depending on the context. Moreover, recall that $e$ and $g_r$ are vector fields on $Q$ that are in one-to-one correspondence with their counterparts on $\T Q$ of \ms, cf. Section \ref{sec:MS}.

In order to formulate the main result we need to define a mechanical analogue of the (vector) relative degree, called also characteristic number, cf. textbooks \cite{nij, Isidori}.

\begin{definition}
	\label{def:reldeg}
	The mechanical control system \ms {equipped} with $\mathbb{R}^m$-valued configuration output map $h:Q\rightarrow \mathbb{R}^m$, $h(x)=\left( h_1(x),\ldots,h_m(x)\right)^T $ has the vector relative half-degree $\left(\nu_1,\ldots,\nu_m \right) $ around $x_0$ if
	\begin{enumerate}[(i)]
		\item $L_{g_r}L_e^{q}h_\ell=0$,\qquad		for $1\leq \ell,r \leq m$ \ and \  {$0\leq q \leq \nu_\ell-2$},
		\item the $m\times m$ decoupling matrix, of the virtual system $\dot{x}=e(x)+\sum_{r=1}^{m}g_r(x)u_r$, $y=h(x)$, given by
		\begin{align*}
			D(x)=\left(L_{g_r}L_e^{\nu_\ell-1} h_\ell \right) (x),
		\end{align*}
		is of full rank equal to $m$, around $x_0$.
	\end{enumerate}
\end{definition}

	\begin{remark}
		Notice that the decoupling matrix $D$ consists of functions $L_{g_r}L_e^{\nu_{\ell}-1}h_\ell$ on $Q$ while the decoupling matrix $\mathbb{D}$, {introduced} in Section \ref{sec:PF}, consists of $L_{G_r}L_F^{\rho_\ell-1}h_\ell$ that are functions on $\T Q$ defined with the help of $G_r=g_r(x)\fp{}{v}$ and $F=v\fp{}{x}-(v^T\Gamma(x)v +e(x))\fp{}{v}$, which are the control vector fields and the drift of \ms, respectively.
\end{remark}

\textbf{Example 1 (cont.):}
	For system \eqref{eq:ex11} of Example 1 we have $e=x^3\fp{}{x^2}$, $g_1=\fp{}{x^1}$, and $g_2=g^2_2\fp{}{x^2}+\fp{}{x^3}$. If $g_2^2\neq0$, then $\nu_1=\nu_2=1$ and we have $\rho_1=2\nu_1$, $\rho_2=2\nu_2$.
	 If $g_2^2\equiv0$, then $\nu_1=1$, $\nu_2=2$, independently of the values of $\Gamma^2_{22}$ and $\Gamma^2_{33}$, since the relative half-degrees depend on $e$ and $g_r$'s only. However, the vector relative degree $(\rho_1,\rho_2)$ is either $(2,3)$ or $(2,4)$ as discussed in Example 1.

In the Appendix, we define the covariant derivative $\nabla \omega$ of a differential 1-form $\omega$ and, in particular, we apply it to $\omega=dL^k_eh_\ell$ and use $\nabla(dL^k_eh_\ell)$ in the following theorem.

\begin{theorem}
	\label{thm:MF-out}
	For the mechanical control system \ms with outputs the MIOLD-problem is solvable, locally around $x_0\in Q$ and globally in $v\in \T_xQ$, if and only if the system \ms
	\begin{enumerate}[(\text{MR}1)]	
		\item has a well defined vector relative half-degree $\left(  \nu_1,\ldots,\nu_m \right)  $,  i.e.
		\begin{align*}
			\rk D(x_0) = m;
		\end{align*}
		\item satisfies
			\begin{align*}
		{	\nabla(dL^{q}_eh_\ell)}&=0\qquad \text{ for } 1\leq \ell \leq m \text{ and } 0\leq q \leq \nu_\ell-2.
		\end{align*}
	\end{enumerate}
\end{theorem}

\begin{remark}
	Condition (MR2) ensures that the differentials of the output functions $h_\ell$ and their successive $\nu_{{\ell}}-2$ Lie derivatives with respect to $e$ are covariantly constant. Equivalently, the new coordinates, given by $h_\ell$ and their successive Lie derivatives $L_e^kh_\ell$, are covariantly linear (precisely, except the last one of each chain of integrators, that can be made covariantly linear by an appropriate feedback). This gives a geometric insight into the mechanical linearization procedure.
\end{remark}
\begin{remark}
	Both conditions (MR1) and (MR2) can be tested using differentiation and algebraic operations only and involve objects defined on $Q$ (of dimension $n$) and not on $\T Q$ (of dimension $2n$). Condition (MR1) uses vector fields $e$ and $g_r$ of the virtual system $\dot{x}=e(x)+\sum_{r=1}^{m}g_r(x)u_r$ only while to calculate the covariant derivatives of condition (MR2) we use also the Christoffel symbols $\Gamma^i_{jk}$.
\end{remark}

\begin{proof}	Necessity. 	 	First, we will show that the conditions hold for the mechanical input-output linear and decupled system \eqref{eq:md}, then we will prove that they are invariant under mechanical diffeomorphisms and feedback \eqref{eq:diff}-\eqref{eq:feed}. 
	
	It is immediate to see that the vector relative half-degree of \eqref{eq:md} is $\left(  \nu_1,\ldots,\nu_m \right)  $ and $D=I_m$, the identity matrix,  so (MR1) holds.
	To avoid summation over double indices now we will rewrite \eqref{eq:md} as follows. Set $\mu_0=0$, $\mu_\ell=\sum_{i=1}^\ell\nu_i$, for $1\leq \ell \leq m$, and denote $\mu=\mu_m$. Then the system \eqref{eq:md} reads (we drop the "tildas"):
\begin{align}
		\label{eq:itnf}
		\begin{split}
		&\dot{x}^i=v^i, \\
		&\dot{v}^i=x^{i+1},   \\
		&\dot{v}^{\mu_\ell}=u_{\ell},\\ 
		&\dot{v}^{i}=-\Gamma^i_{jk}(x)v^jv^k+e^i(x)+\sum_{r=1}^{m}g^i_r(x) u_r, \\
		&y_\ell=h_\ell(x)=x^{\mu_{\ell-1}+1},
		\end{split}
	\begin{split}
		1\leq & \ i \leq n, \\
		\mu_{\ell-1}+1\leq & \ i \leq \mu_{\ell}-1,\\
		1\leq & \ \ell \leq m,\\
		 \mu+1\leq & \  i\leq n \\
		 1\leq & \ \ell \leq m  .
	\end{split}
\end{align}

\noindent Now observe that in \eqref{eq:itnf} we have  $L_e^{l-1}h_\ell=x^{\mu_{\ell-1}+l}$, $1\leq l \leq \nu_\ell$, and that the Christoffel symbols in the equations for $\dot v^i$, where $1\leq i \leq \mu$, vanish identically. Setting $q=l-1$ we conclude, {for $0\leq q \leq \nu_{{\ell}}-2$, that}
\begin{align}
	\label{eq:call}
	\nabla( dL_e^{q}h_\ell)=\nabla( dL_e^{l-1}h_\ell)= {\fp{^2}{x^j\partial x^k}}(x^{\mu_{\ell-1}+l})\ dx^j\otimes dx^k= 0, \quad \text{for} \quad 1\leq j, k\leq n,
\end{align} 
 see \eqref{eq:cdwa} in the Appendix, so (MR2) holds as well.
	
	The invariance under diffeomorphisms is obvious since both the Lie derivative and the covariant derivative are  geometric operations and the notion of relative degree does not depend on coordinates. Now, we will prove invariance under the mechanical feedback \eqref{eq:feed}. Note that the objects defining \ms, namely $\Gamma^i_{jk}$, $e$, and $g_r$, change under feedback in the following way
	\begin{align}
		\label{eq:Ch}
		\tilde{\Gamma}^i_{jk}=\Gamma^i_{jk}-\sum_{r=1}^m g_r^i\gamma^r_{jk}, \quad
		\tilde{e}=e+\sum_{r=1}^m g_r\alpha^r, \quad
		\tilde{g}_s=\sum_{r=1}^m \beta_s^r g_r,
	\end{align}
	and define the feedback modified system $(\widetilde{\mathcal{MSO}})$ with the maintained output $y=h(x)$. The decoupling matrix $\widetilde D( x)$, of the virtual system $\dot{ x}=\ti e( x) + \sum_{r=1}^{m} \ti g_r( x)\ti u_r $ of $(\widetilde{ \mathcal{MSO}})$, is $\widetilde D( x)=\left( L_{\ti g_r}L_{\ti e}^{\ti \nu_{\ell}-1 }h_\ell\right) ( x)$.
	It is known, see e.g. \cite{nij,Isidori}, that for $0\leq q \leq \nu_\ell-1$
	\begin{align}
		\label{eq:gw}
		L_{\ti e}^qh_\ell=L_e^qh_\ell \quad \text{and} \quad L_{\ti g_s}L_{\ti e}^{\nu_\ell-1}h_\ell=\sum_{r=1}^{m}\beta^r_s L_{g_r}L_e^{\nu_\ell-1}h_\ell.
	\end{align}
Item (i) of Definition~\ref{def:reldeg} implies that the former equality of \eqref{eq:gw} yields
 $L_{\ti g_r}L_{\ti e}^{q}h_\ell=0$, for $0\leq q\leq \nu_\ell-2$, and item (ii) of Definition~\ref{def:reldeg} implies that the latter yields $\widetilde D(x): =\left( \beta^r_s L_{g_r}L_e^{\nu_\ell-1}h_\ell\right) =D(x)\beta(x)$. Therefore (i) and (ii) are feedback invariant showing that (MR1) is necessary.
	To prove necessity of (MR2), by a direct calculation we show that the covariant derivative $\ti \nabla$, whose Christoffel symbols $\ti \Gamma^i_{jk}$ are given by \eqref{eq:Ch}, acts on a differential one-form $\omega\in \Lambda(Q)$ as follows
	\begin{align*}
		\ti{\nabla}\omega=\nabla\omega + \sum_{r=1}^{m }\omega(g_r) \cdot \gamma^r,
	\end{align*}
	where the matrices $\gamma^r=\left( \gamma^r_{jk}\right)$, for $1\leq j,k \leq n$. Thus, for $\omega=dL^q_{\ti{e}}h_\ell$, by \eqref{eq:gw} we have
	\begin{align*}
		\ti{\nabla}dL^q_{\ti{e}}h_\ell=\ti{\nabla}dL^q_{e}h_\ell= \nabla dL^q_{e}h_\ell+ \sum_{r=1}^{m} dL_e^qh_\ell (g_r)\gamma^r=0,
	\end{align*}
	for $0\leq q \leq \nu_\ell-2$,	where the first term vanishes by \eqref{eq:call} and the second by the definition of the relative half-degree since $dL_e^qh_\ell (g_r)=L_{g_r}L_e^qh_\ell$.

	Sufficiency.
	Consider a mechanical system with outputs \ms satisfying (MR1)-(MR2). We will transform it to the form~\eqref{eq:itnf}, which is an equivalent way of representing \eqref{eq:md}. By (MR1), the mechanical system \ms has a well defined vector relative half-degree $\left(  \nu_1,\ldots,\nu_m \right)  $, so set $\mu_0,\mu_1,\ldots,\mu_m=\mu$ as at the beginning of the necessity part and define functions $\phi^1,\ldots,\phi^\mu$ by {taking}, for $\ell=1,\ldots,m$,
		\begin{align*}
			\phi^{\mu_{\ell-1}+l}=L_e^{l-1}h_\ell, \qquad \qquad 1\leq l \leq \nu_\ell.
	\end{align*}
	It is well known that the functions $\phi^1,\ldots,\phi^\mu$ are locally independent around $x_0$, see e.g. \cite{nij,Isidori}, and complete them by $\phi^{\mu+1},\ldots,\phi^n$ so that $\phi=(\phi^1,\ldots,\phi^n)$ is a local diffeomorphism.
	 Set $\ti x^i=\phi^i(x)$ and $\ti v^i= \fp{\phi^i}{x}\cdot v$, for $1\leq i \leq n$. In $(\ti x,\ti v)$-coordinates the system \ms reads
	\begin{align*}
	\begin{array}{lc}
		\dot{\ti x}^i=\ti v^i, & \qquad \quad 1\leq  \ i \leq n,\\
		\dot{\ti v}^i=-\ti \Gamma^i_{jk} \ti v^j \ti v^k+ \ti x^{i+1},   & \qquad \mu_{\ell-1}+1\leq  \ i \leq \mu_{\ell}-1,\\
		\!\!\!	\dot{\ti v}^{\mu_\ell}=-\ti \Gamma^{\mu_\ell}_{jk} \ti v^j \ti v^k+ \ti e^{\mu_\ell} + \sum\limits_{r=1}^{m}\ti d^{\mu_\ell}_r u_r, & \qquad \quad 1\leq   \ell \leq m,\\ 
		\dot{\ti v}^{i}=-\ti \Gamma^i_{jk}\ti v^j\ti v^k+\ti e^i+\sum\limits_{r=1}^{m}\ti g^i_r u_r, & \quad  \mu+1\leq  \  i\leq n \\
		y_\ell=\ti h_\ell(\ti x)=\ti x^{\mu_{\ell-1}+1}, & \qquad \quad	1\leq   \ell \leq m
	\end{array}
\end{align*}
where $\ti \Gamma^i_{jk}$ are Christoffel symbols expressed in $\ti x$-coordinates and $\left(\ti d^{\mu_\ell}_r \right)=\ti D $ is the invertible, due to (MR1), decoupling matrix $D$ expressed in $\ti x$-coordinates, that is, $\ti D(\ti x)=D(\phi^{-1}(\ti x))$. We will prove that $\ti \Gamma^i_{jk}=0$, for $1\leq j,k\leq n$ and all $1\leq i\leq n$ such that $i\neq \mu_\ell $, for $1\leq \ell \leq m$; any such index $i$ is of the form $i=\mu_{\ell-1}+l$ for a certain $1\leq l \leq \nu_{\ell}-1$. To this end, calculate in $\ti x$-coordinates using (MR2)
\begin{align*}
	0=\nabla\left(dL_e^{l-1}h_{\ell} \right)= \nabla\left(d\phi^{\mu_{\ell-1}+l} \right)= \nabla\left(d\ti x^{\mu_{\ell-1}+l} \right)
\end{align*}
and denoting $i=\mu_{\ell-1}+l$ and using \eqref{eq:cdwa} we get 
\begin{align*}
	0=\nabla\left(d\ti x^i \right)= \fp{^2\ti x^i}{\ti x^j\partial \ti x^k}+\ti \Gamma^s_{jk} \fp{\ti x^i}{\ti x^s}=0+\ti \Gamma^s_{jk} \delta^i_s=\ti \Gamma^i_{jk} 
\end{align*}
implying that $\ti \Gamma^i_{jk}=0$ and thus giving $\dot{\ti v}^i=\ti x^{i+1}$, for any $1\leq i \leq n$, $i\neq \mu_{\ell} $. Applying the invertible feedback 
\begin{align*}
	\ti u_\ell=  \ti v^T\ti \gamma^\ell(\ti x) \ti v+\ti \alpha^\ell(\ti x)+\sum_{r=1}^{m}\ti \beta^\ell_r(\ti x)u_r,
\end{align*}
where $\ti v^T\ti \gamma^\ell(\ti x) \ti v= -\ti \Gamma^{\mu_\ell}_{jk} \ti v^j \ti v^k $, $\ti \alpha^\ell=\ti e^{\mu_\ell}$, and $\ti \beta^\ell_r= \ti d^{\mu_\ell}_r$ (implying that the feedback is, indeed, invertible since so is $\ti D=\left(\ti d^{\mu_\ell}_r \right)$), we get $\dot{\ti v}^{\mu_\ell}=\ti u_\ell$, for $1 \leq \ell \leq m$, and thus we obtain the desired form \eqref{eq:itnf} in $(\ti x,\ti v)$-coordinates and with respect to the controls $\ti u=(\ti u_1,\ldots,\ti u_m)^T$.

\end{proof}

It is worth to interpret the linearizability conditions. The first condition (MR1) ensures that in the $(\ti{x},\ti v)$-coordinates the control appears in  equations for $\dot{	\ti  v}^{\nu_\ell}_\ell$ and $\dot{	\ti  v}_{m+1}$ only, see \eqref{eq:md}. The second condition (MR2) ensures that by introducing the new coordinates $(\ti x,\ti v)$  we compensate the Christoffel symbols in the $(\ti x_o,\ti v_o)$-subsystem.

To summarize, the mechanical diffeomorphism and feedback that perform input-output linearization and non-interacting are given by
\begin{align}
	\label{eq:MF-diff}
	&\left(\ti{x},\ti{v} \right) =\Phi(x,v)=\left(\phi(x),\fp{\phi}{x}(x)v\right) \quad \text{with} \\
	&\phi(x)=\big(h_1,L_eh_1,\ldots,L^{\nu_1-1}_eh_1,h_2,L_eh_2,\ldots,L^{\nu_2-1}_eh_2,  \ldots, h_m,L_eh_m,\ldots,L^{\nu_m-1}_eh_m, \phi^{\mu+1}, \ldots,\phi^n \big)^T\nonumber,
\end{align}
and the feedback modified controls, for $1\leq \ell \leq m,$ are
\begin{align*}
	\tilde{u}_\ell=\left(\fp{^2L_e^{\nu_i-1}h_\ell}{x^j\partial x^k}-\fp{L_e^{\nu_i-1}h_\ell}{x^i}\Gamma^i_{jk} \right)v^jv^k+L^{\nu_\ell}_eh_\ell+\sum_{r=1}^m L_{g_r}L^{\nu_\ell-1}_eh_\ell u_r,
\end{align*}
or using matrix notation
\begin{align*}
	\tilde{u}=\mathcal{C}(x,v)+\mathcal{A}(x)+D(x) u,
\end{align*}
where $u=\left(u_1,\ldots, u_m \right)^T $, $\tilde{u}=\left(\tilde{u}_1,\ldots, \tilde{u}_m \right)^T $, 
\begin{align*}
	\mathcal{A}(x)=\left(\begin{array}{c}
		L^{\nu_1}_eh_1\\ 
		\vdots\\ 
		L^{\nu_m}_eh_m
	\end{array}  \right), \qquad 
	D(x)=\left(\begin{array}{ccc}
		L_{g_1}L^{\nu_1-1}_eh_1&\ldots&L_{g_m}L^{\nu_1-1}_eh_1\\ 
		\vdots& \ddots &\vdots\\ 
		L_{g_1}L^{\nu_m-1}_eh_m&\ldots&L_{g_m}L^{\nu_m-1}_eh_m
	\end{array}  \right),
\end{align*}
\begin{align*}
	\mathcal{C}(x,v)=\left(\begin{array}{c}
		\left(\fp{^2L_e^{\nu_1-1}h_1}{x^j\partial x^k}-\fp{L_e^{\nu_1-1}h_1}{x^i}\Gamma^i_{jk} \right)v^jv^k\\ 
		\vdots\\ 
		\left(\fp{^2L_e^{\nu_m-1}h_m}{x^j\partial x^k}-\fp{L_e^{\nu_m-1}h_m}{x^i}\Gamma^i_{jk} \right)v^jv^k
	\end{array}  \right),
\end{align*}
hence the MIOLD-linearizing mechanical feedback reads
\begin{align}
	\label{eq:MF-u}
	u=D^{-1}(x) \left(-\mathcal{C}(x,v)-\mathcal{A}(x)+\tilde{u} \right) .
\end{align}

\begin{corollary}
	If the control distribution $\mathcal{G}=\spn{g_1,\ldots,g_m}$ is involutive {around} $x_0$ (for instance, if the system has a single control), then it is possible to choose the functions $\phi^{\mu+1}(x), \ldots,\phi^n(x)$ such that
	\begin{align*}
		L_{g_r}\phi^i=0, \qquad \qquad \text{for } 1\leq r \leq m, \ \mu+1\leq i \leq n,
	\end{align*}
	and in coordinates $(\ti x^1,\ldots,\ti x^\mu,\ti x^{\mu+1},\ldots,\ti x^n)$, where $\ti x^i=\phi^i$, $\mu+1\leq i \leq n$ (completed {by} the corresponding velocities $\ti v^i=\dot{\ti x}^i$), the unobserved part of \eqref{eq:md} does not contain controls and takes the form
	\begin{align*}
		&	\dot{	\ti  x}_{m+1}=\ti v_{m+1}\\
		&	\dot{	\ti  v}_{m+1}= -\ti v^T \ti \Gamma(\ti x)\ti v + \ti e(\ti x) .
	\end{align*}
\end{corollary}
The proof is a direct application of the Frobenius theorem and follows the same line as that of  Proposition 5.1.2. in \cite{Isidori}.

{If the unobserved $(\ti x_{m+1}, \ti v_{m+1})$-part is absent in \eqref{eq:md}, then the system takes the following form.}

\begin{corollary}
	\label{coro:MF-out}
	The mechanical control system \ms is, locally around $x_0\in Q$ and globally in $v$, equivalent via a mechanical feedback transformation of the form \eqref{eq:diff}-\eqref{eq:feed} to \eqref{eq:md} with $(\ti x_{m+1}, \ti v_{m+1})$-part absent if and only if \ms satisfies conditions (MR1)-(MR2) of Theorem \ref{thm:MF-out} and, moreover,
	\begin{align*}
		\mu=\sum_{\ell=1}^m \nu_\ell=n.
	\end{align*}
{In this case, the system is linear on the whole state-space and consist of $m$ chains of $2\nu_{\ell}$-fold integrators (which is the canonical form of any controllable linear mechanical control system, see \cite{lms}) with outputs equal to the top configuration variables of each chain.} 
\end{corollary}


{
Now, for mechanical system \eqref{eq:ms_bez} without outputs
\begin{align}
	\begin{split}
		\dot{x}&=v\\
		\dot{v}&=-v^T\Gamma(x)v+e(x)+\sum_{r=1}^{m}g_r(x) u_r,
	\end{split}
	\label{eq:ms_bez}
	\tag{$\mathcal{MS}$}
\end{align}
consider the problem of} mechanical feedback linearization, shortly MF-linearization \cite{aut,tac}, i.e. the problem of whether we can transform the system \eqref{eq:ms_bez} without output to a linear mechanical system using transformations \eqref{eq:diff}-\eqref{eq:feed}.

\begin{proposition}
	\label{prop:MF-lin}
	The mechanical control system \eqref{eq:ms_bez} is, locally around $x_0\in Q$, MF-linearizable if and only if there exist $m$ functions $h_1,\ldots,h_m\in C^{\infty}(Q)$ for which
	\begin{enumerate}[(MF1)]	
		\item the vector relative half-degree $\left(  {\nu}_1,\ldots,{\nu}_m \right)  $ of \eqref{eq:ms_bez}, with $h=(h_1,\ldots,h_m)$ considered as its output, is well defined and satisfies $\sum\limits^m_{\ell=1}{\nu}_\ell=n$.
		\item $		\nabla(dL^k_eh_\ell)=0$ \qquad for  $0\leq k \leq {\nu}_\ell-2$.
	\end{enumerate}
\end{proposition}

\begin{remark}
The general nonlinear control system $\dot z= F(z)+\sum_{r=1}^{m}G_r(z)u_r$, $z\in Z\subset \mathbb{R}^N$, is locally feedback linearizable if and only if there exist $m$ functions $h_1,\ldots,h_m$ whose vector relative degree $(\rho_1,\ldots,\rho_m)$ is well defined and satisfies $\sum_{\ell=1}^{m}\rho_\ell=N$, cf. Lemma~5.2.1 in textbook \cite{Isidori}. Proposition \ref{prop:MF-lin} can thus be considered as a mechanical counterpart of that result and conditions (MF1)-(MF2) guarantee that the linearizing transformations exist and are mechanical.
\end{remark}

It is well known that static feedback linearizable systems $\dot{z}=F(z,u)$ are flat \cite{fli,lev,martin,murr}, that is, we can find $\varphi=(\varphi_1,\ldots,\varphi_m)$ called a flat output, where $\varphi_\ell=\varphi_\ell(z,u,\ldots,u^{(\kappa)})$, such that
\begin{align}
	\label{eq:flat}
z=\zeta(\varphi,\dot \varphi, \ldots, \varphi^{(\sigma-1)}) \quad \text{and} \quad u=\delta(\varphi,\dot \varphi, \ldots, \varphi^{(\sigma)}),
\end{align}
for some integers $k\geq -1$ and $\sigma\geq 1$.
 If the mechanical system \eqref{eq:ms_bez} without output, is MF-linearizable then is thus flat  and the flatness property reflects additionally its mechanical structure. Namely, as a flat output of \eqref{eq:ms_bez} we can take $(\varphi_1,\ldots,\varphi_m)=(h_1,\ldots,h_m)$ given by Proposition \ref{prop:MF-lin} (that holds for \eqref{eq:ms_bez} since it is MF-linearizable). Moreover, expressing the configurations $x^i$ involves time-derivatives of $h_\ell$ of even orders only and that of velocities $v^i$ {of linear combination of even and odd orders}. The differential weight \cite{nic1,nic2}, which is the number time-derivatives $\varphi_\ell^{(j)}$ involved in \eqref{eq:flat}, of the flat output $h=(h_1,\ldots,h_m)$ is $2n+m$, which is the minimal possible thus confirming that the system is static feedback linearizable (cf. Theorem 2.2 in \cite{nic2}).   
More precisely, we have the following result. Given local coordinates $(x^i)$ on $Q$, the induced local coordinates on $\T Q$ defined by $(x^i,\dot x^i)=(x^i,v^i)$ will be called mechanical state-variables of \eqref{eq:ms_bez} and consists of configuration variables $x=(x^i)$ and corresponding velocities $v=(v^i)=(\dot x^i)$.

\begin{proposition}
	\label{prop:flat}
	If the mechanical control system \eqref{eq:ms_bez} is, locally around $x_0$, MF-linearizable
	then $h=(h_1,\ldots,h_m)$ of Proposition \ref{prop:MF-lin} is a flat output of differential weight $2n+m$ and, moreover, \eqref{eq:ms_bez} is locally MF-equivalent to \eqref{eq:md}, with  $(\ti x_{m+1}, \ti v_{m+1})$-part absent, where
	\begin{align*}
		\ti x&=\left( h_1, \ddot{h}_1,\ldots, h^{(2\nu_1-2)}_1, h_2, \ddot{h}_2 \ldots, h_m^{(2\nu_m-2)}\right),\\
	\ti	v&={\left( \dot{ h}_1, \dddot{h}_1,\ldots, h^{(2\nu_1-1)}_1,  \dot{h}_2, \dddot h_2,\ldots, h_m^{(2\nu_m-1)}\right)},
	\end{align*}
	where $\left(\nu_1,\ldots,{\nu_m}\right)$ is the vector relative half-degree of \eqref{eq:ms_bez} with $h=(h_1,\ldots,h_m)$ considered as output.
Any other mechanical state-variables $(x,v)$ of \eqref{eq:ms_bez} are given by
	\begin{align*}
	 x&=\zeta\left( h_1, \ddot{h}_1,\ldots, h^{(2\nu_1-2)}_1, h_2, \ddot{h}_2 \ldots, h_m^{(2\nu_m-2)}\right),\\
	v&=\varrho{\left(  h_1, \dot{h}_1,\ldots, h^{(2\nu_1-1)}_1,  h_2, \dot h_2,\ldots, h_m^{(2\nu_m-1)}\right)},
\end{align*}
where the map $(\zeta,\varrho)$ is locally invertible and $\varrho$ is a linear map of odd derivatives $h^{(2j-1)}_\ell$ with coefficients depending on even derivatives $h^{(2j)}_\ell$ only.
\end{proposition}
A proof is a direct consequence of the sufficiency part of the proof of Theorem \ref{thm:MF-out}.

\begin{remark}
Any configuration variables $x=(x^i)$ are expressed in terms of even derivatives $h^{(2j)}_\ell$ only. On the other hand, the velocity variables $\ti v=(\ti v^i)$ of \eqref{eq:md} are expressed in terms of odd derivatives  $h^{(2j-1)}_\ell$ only, while others velocity variables  $v=(v^i)$ as their linear functions (with coefficients depending on even derivatives). In fact, $v=\sum\limits_{\ell=1}^{m}\sum\limits_{j=0}^{\nu_{{\ell}}-1}\fp{\zeta\ \ }{h_\ell^{(2j)}}h_\ell^{(2j+1)}=$ $= \sum\limits_{\ell=1}^{m}\sum\limits_{j=0}^{\nu_{{\ell}}-1}\fp{\zeta\ \ }{h_\ell^{(2j)}}\ti v_\ell^{j+1}$ and, indeed, the coefficients $\fp{\zeta\ \ }{h_\ell^{(2j)}}$ depend on $h_\ell^{(2j)}$ only.
\end{remark}

\section{Examples}
\subsection{The Inertia Wheel Pendulum}
Consider the equation of dynamics of the Inertia Wheel Pendulum \cite{iwp}, i.e. a pendulum with a rotating wheel attached, as depicted in Figure \ref{fig:iwp},
\begin{align*}
	\begin{split}
		\dot{x}^1&=v^1\\
		\dot{x}^2&=v^2
	\end{split}
	\begin{split}
		\dot{v}^1&=e^1+g^1 u\\
		\dot{v}^2&=e^2+g^2 u,
	\end{split}
\end{align*}
 where $	e^1=\frac{m_0}{m_d}\sin x^1$,  $g^1=-\frac{1}{m_d}$, $e^2=-\frac{m_0}{m_d}\sin x^1$,  $g^2=\frac{m_d+J_2}{J_2 m_d}$, and $m_0=aL_1(m_1+2m_2)$, $m_d=L_1^2(m_1+4m_2)+J_1$, $J_2$ are constant parameters. Above, {$x^1=\theta^1\in \mathbb{S}^1$} denotes the angle of the pendulum measured from the vertical position, and {$x^2=\theta^2\in \mathbb{S}^1$} is the angle of the wheel. The masses and the momenta of inertia of the pendulum and the wheel are $m_1, \ m_2$ and $J_1, \ J_2$, respectively. The distance to the center of the pendulum is $L_1$ and $a$ is the gravitational constant. The  torque applied to the wheel is the control signal. 
\begin{figure}[h!]
	\centering
	\includegraphics[width=0.2\linewidth]{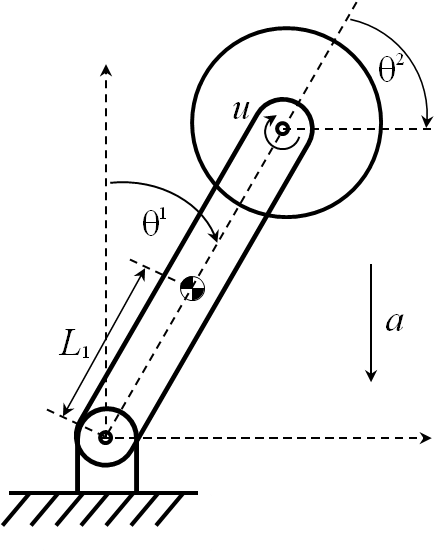}
	\caption{The Inertia Wheel Pendulum}
	\label{fig:iwp}
\end{figure}
In the following calculation we use vector fields on {$Q=\mathbb{S}^1\times\mathbb{S}^1$}, namely $e(x)=e^i(x)\fp{}{x^i}$ and $g(x)=g^i(x)\fp{}{x^i}$ with $e^i(x)$ and $g^i(x)$ defined above.
The output of the system is $y=h(x)=x^1$. It is straightforward to see that the system has the relative half-degree $\nu=1$, since $L_gh=-\frac{1}{m_d}\neq 0$. Therefore, by applying the mechanical feedback $\ti u=e^1+g^1 u$ we get the normal form \eqref{eq:md}:
\begin{align*}
	\begin{split}
		y=x^1\\
		\\
		\\
	\end{split}
	\begin{split}
		\dot{x}^1&=v^1\\
		\\
		\dot{x}^2&=v^2
	\end{split}
\qquad \qquad
	\begin{split}
		\dot{v}^1&=e^1+g^1 \frac{1}{g^1}(-e^1+\ti{u})=\ti{u}\\
		\dot{v}^2&=e^2+g^2 \frac{1}{g^1}(-e^1+\ti{u})=e^2-\frac{g^2}{g^1}e^1+\frac{g^2}{g^1}\ti u=\ti e^2+ \ti g^2 \ti u.
	\end{split}
\end{align*}
The $2$-dimensional subsystem $(x^1,v^1)$, with the output $y=x^1={\theta^1}$, is linear and mechanical and decoupled from the unobservable $(x^2,v^2)$-part that is still mechanical. This directly corresponds to Theorem \ref{thm:MF-out} (note that for $\nu=1$ condition (MR2) is empty).

To illustrate condition (MR2), consider the same Inertia Wheel Pendulum with the output $y=h(x)=\frac{m_d+J_2}{J_2}x^1+x^2$. We have the relative half-degree $\nu=2=n$, since $L_gh=0$ and $L_gL_eh\neq 0$. Moreover, $\nabla dh=\fp{^2 h}{x^j\partial x^k}-\Gamma^k_{ji}\fp{h}{x^k}=0$, since $\fp{^2 h}{x^j\partial x^k}=0$ and all Christoffel symbols $\Gamma^k_{ji}$ are zero. Thus the system satisfies (MR1)-(MR2) of Theorem \ref{thm:MF-out} and is locally equivalent (around any $x_0$ satisfying $x^1_0\neq\pm\frac{\pi}{2}$), via a mechanical feedback transformation of the form \eqref{eq:diff}-\eqref{eq:feed}, to the system \eqref{eq:md} without the unobserved part.

If we consider the MF-linearization problem of the Inertia Wheel Pendulum, then obviously the above function $h(x)=\frac{m_d+J_2}{J_2}x^1+x^2$ satisfies conditions (MF1)-(MF2) of Proposition \ref{prop:MF-lin} (all others being of the form $ch(x)$, where $c\in \mathbb{R}^*$).
Indeed, the mechanical diffeomorphism  $(\ti{x},\ti{v})=\Phi(x,v)=(\phi(x),\fp{\phi}{x}(x) v)$ with $\phi(x)=(h,L_eh)^T$ MF-linearizes the system. The system in new coordinates reads
\begin{align}
	\label{eq:iwp-u}
	\dot{\tilde{x}}^1&=\frac{m_d+J_2}{J_2}v^1+v^2=\tilde{v}^1\nonumber\\
	\dot{\tilde{v}}^1&=\frac{m_d+J_2}{J_2}\left(\frac{m_0}{m_d}\sin x^1 - \frac{1}{m_d}u \right) -\frac{m_0}{m_d} \sin x^1 + \frac{m_d +J_2}{m_2 J_2}u\nonumber\\
	&=\frac{m_0}{J_2}\sin x^1=L_eh=\tilde{x}^2\\
	\dot{\tilde{x}}^2&=\frac{m_0}{J_2}\cos x^1 v^1=\tilde{v}^2\nonumber\\
	\dot{\tilde{v}}^2&=-\frac{m_0}{J_2} \sin x^1 v^1v^1+\frac{m_0^2}{2m_dJ_2}\sin(2x^1)- \frac{m_0}{m_d J_2}\cos x^1u=\tilde{u}\nonumber,
\end{align}
which is linear, mechanical, and consists of a chain of 4 integrators, {assuming that $x^1\neq \pm \frac{\pi}{2}$}. Finally, consider another output given by a nonlinear function of the previous one, namely $\ti h=\ti h\left( \frac{m_d+J_2}{J_2}x^1+x^2\right)$, {where $\ti h$ is not $ch$, $c\in \mathbb{R}^*$}. Now the system is input-output decouplable (in fact, feedback linearizable), since the relative degree $\rho=4$, however not mechanically input-output decouplable. Indeed, the relative half-degree is still $\nu=2$, yet $\nabla d\ti h\neq 0$ which contradicts (MR2) of Theorem \ref{thm:MF-out}.   

\subsection{TORA3}
We will study {the MIOLD-problem and its relation with MF-linearization} of the TORA3 system (see Figure \ref{fig:tora}) proposed in \cite{tac}, which is based on the TORA system (Translational Oscillator with Rotational Actuator) studied in the literature, e.g. \cite{tora} (however we add gravitational effects). It consists of a two dimensional spring-mass system, with masses $m_1, m_2$ and spring constants $k_1, k_2$, respectively. A pendulum of length $l_3$, mass $m_3$, and moment of inertia $J_3$ is added to the second body. The displacements of the bodies are denoted by $x^1$ and $x^2$, respectively, and the angle of the pendulum by $x^3$. The control $u$ is a torque applied to the pendulum. The equations of the mechanical system  on $\T Q$, {where $Q=\mathbb{R}^2\times \mathbb{S}^1$,} read
\begin{align*}
	\begin{split}
		&\dot{x}^1=v^1\qquad \dot{v}^1=e^1\\
		&\dot{x}^2=v^2\qquad \dot{v}^2=-\Gamma^2_{33}v^3 v^3+e^2+g^2 u\\
		&\dot{x}^3=v^3\qquad \dot{v}^3=-\Gamma^3_{33}v^3 v^3+e^3+g^3 u
	\end{split}
\end{align*}
where ${\Gamma^2_{33}=\frac{-p_0 \sin x^3}{p_1+p_2 \sin^2 x^3}}$, ${\Gamma^3_{33}=\frac{ p_2 \sin x^3 \cos x^3}{p_1+p_2 \sin^2 x^3}}$, ${e^1=-\frac{k_1}{m_1}x^1+\frac{k_2}{m_1} \left(x^2-x^1 \right)}$, ${e^2=\frac{\frac{1}{2}p_2 a \sin 2x^3-p_3(x^2-x^1)}{p_1+p_2 \sin^2 x^3}}$, $e^3=\frac{p_4 \left(x^2-x^1 \right)\cos x^3-p_5 \sin x^3 }{p_1+p_2 \sin^2 x^3}$, ${ g^2=\frac{-m_3 l_3 \cos x^3}{p_1+p_2 \sin^2 x^3}}$, ${g^3=\frac{m_2+m_3}{p_1+p_2 \sin^2 x^3}}$, with constant parameters: 
$p_0=m_3 l_3 (m_3l_3^2+J_3 ), \ p_1=m_2m_3l_3^2+J_3(m_2+m_3)$,
$p_2=m_3^2l_3^2,\\
p_3=k_2\left( m_3l_3^2+J_3\right), \
p_4=m_3l_3 k_2$
$p_5=m_3 l_3 a (m_2+m_3)$. Recall that $e(x)=e^i(x)\fp{}{x^i}$ and $g(x)=g^i(x)\fp{}{x^i}$, with $e^i(x)$ and $g^i(x)$ defined above, are vector fields on $Q$.
\begin{figure}[h!]
	\centering
	\includegraphics[width=0.5\linewidth]{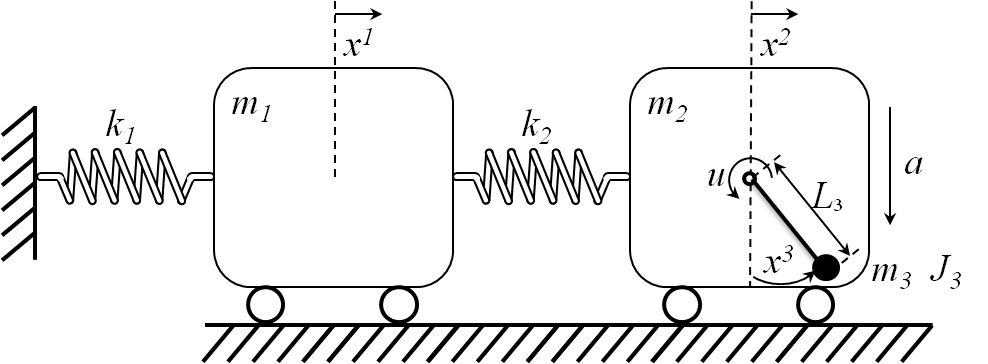}
	\caption{The TORA3 system.}
	\label{fig:tora}
\end{figure}

Set the output function $y=h(x)=x^1$. A straightforward calculation shows that $L_gh=0$, $L_eh=-\frac{k_1}{m_1}x^1+\frac{k_2}{m_3} \left(x^2-x^1 \right)$, and $L_gL_eh\neq0$, so the relative half-degree is $\nu =2$ and, moreover, $\nabla dh=0$. Therefore, applying a change of coordinates $(\ti{x},\ti{v})=\Phi(x,v)=(\phi(x),D\phi(x) v)$ with $\phi(x)=(h,L_eh,x^3)^T$ and the feedback $\ti u=-(\fp{L_eh}{x^2}\Gamma^2_{33}+\fp{L_eh}{x^3}\Gamma^3_{33})v^3v^3+L^2_eh+L_gL_ehu$, the system can be transformed into the following  system of the form \eqref{eq:md}, with $4$-dimensional linear subsystem and $2$-dimensional unobserved part $(\ti x^3,\ti v^3)$,
\begin{align*}
		\begin{split}
		&	y=\ti x^1\\		
	\end{split}
	\begin{split}
		&	\dot{	\ti  x}^1=\ti v^1\\
		&	\dot{	\ti  v}^1=\ti x^2\\
				&	\dot{	\ti  x}^2=\ti v^2\\
		&	\dot{	\ti  v}^2 =\ti u
	\end{split}
	\qquad \qquad \qquad
	\begin{split}
		&	\dot{	\ti  x}^3=\ti v^3\\
		&	\dot{	\ti  v}^3= -\ti \Gamma^3_{33}\ti v^3 \ti v^3+\ti e^3+\ti g^3 \ti u.
	\end{split}
\end{align*}

Choosing the output function as $y=h(x)=\frac{m_1}{m_2+m_3} x^1 +x^2 +\frac{m_3 l_3}{m_2+m_3} \sin x^3$, it is possible to fully MF-linearize the system (see Corollary \ref{coro:MF-out}). A direct calculation shows that now we have $L_gh=L_gL_eh=0$, since $L_eh=-\frac{k_1}{m_2+m_3}x^1$, and $L_gL_e^2h=\frac{k_1 p_4 \cos x^3}{m_1(m_2+m_3)(p_1+p_2 \sin^2x^3)}\neq 0$, {around $x^3_0\neq \pm \frac{\pi}{2}$,} therefore the relative half-degree is $\nu=n=3$. Moreover, $\nabla dh=\nabla dL_eh=0$ thus the system is (fully, input-state) MF-linearizable using a linearizing diffeomorphism $\left( \ti{x},\ti{v}\right)=\left(\phi(x),\fp{\phi}{x}(x)v \right)  $, with $\phi(x)=\left(h,L_eh,L_e^2h \right)^T$ and the feedback is $\ti u=(\fp{^2L_e^2h}{x^jx^k}-\fp{L_e^2h}{x^i}\Gamma^i_{jk})v^jv^k+L^3_eh+L_gL_e^2hu$. The resultant system is the chain of 6 integrators.

\subsection{Forced double pendulum on an oscillating base.}

Consider a forced double pendulum on an oscillating base as depicted on Figure~\ref{fig:2pOb}, see \cite{aut}.
\begin{figure}[h!]
	\centering
	\includegraphics[width=0.5\linewidth]{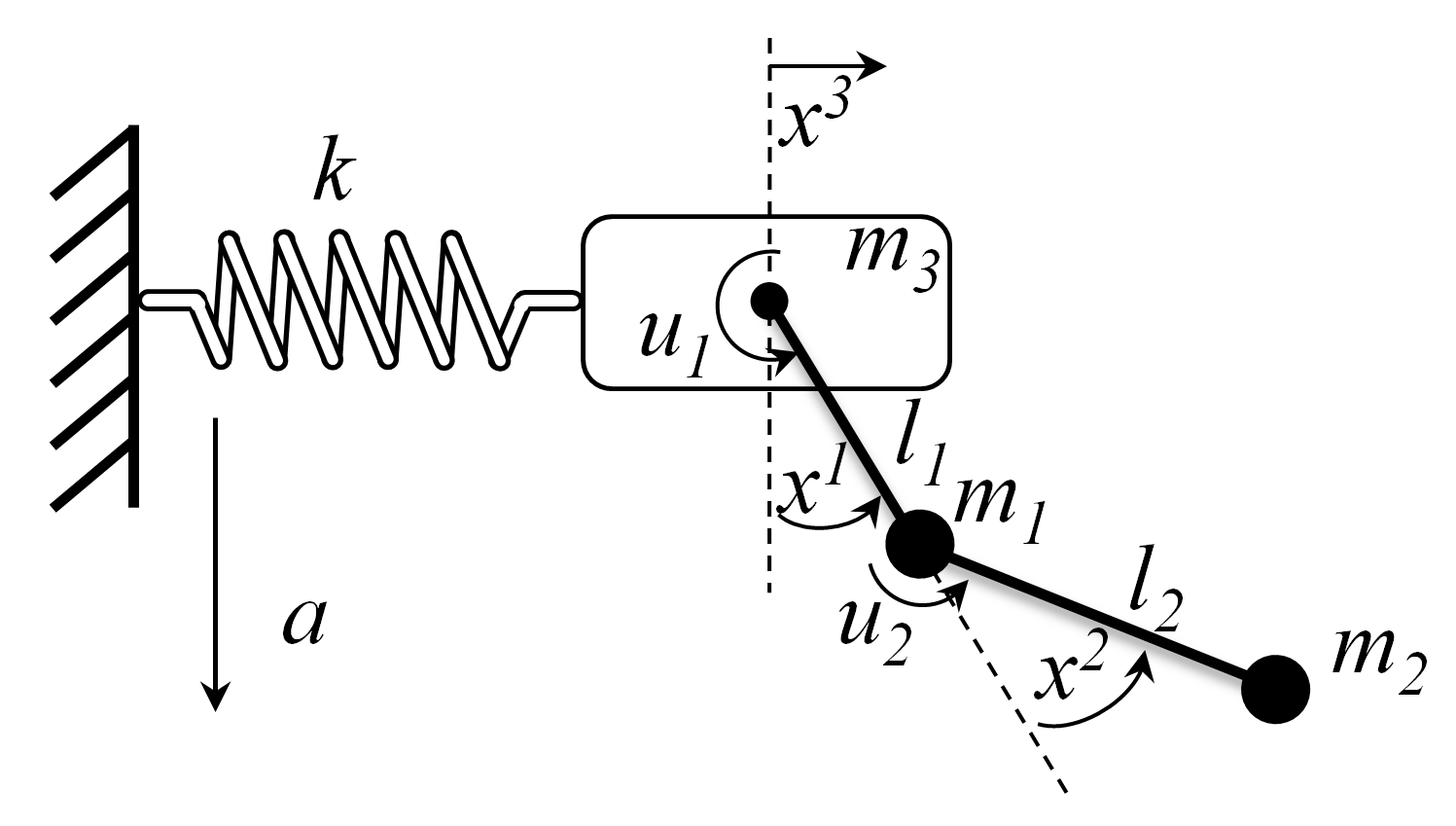}
	\caption{Forced double pendulum on an oscillating base}
	\label{fig:2pOb}
\end{figure}
It consists of a 2-DOF pendulum with actuators in both joints, with masses $m_1$, $m_2$ and lengths $l_1$, $l_2$, respectively. The pendulum is mounted on an harmonic oscillator, with mass $m_3$ and spring constant $k$, that moves horizontally. The displacement of the base is denoted by $x^3$, and the angles of the pendulum by $x^1$ and $x^2$. The controls $u_1$ and $u_2$ are torques applied to, respectively, the first and second joint of the pendulum. 

The kinetic energy is $T=\frac{1}{2}{\rm m}_{ij}\dot{x}^i\dot{x}^j$, where $\rm m_{ij}\geq0$ are elements of the inertia matrix, i.e. the metric (0,2)-tensor, given by:

	\begin{align*}
		&	{\rm m}=\begin{pmatrix}
			l_1^2(m_1+m_2) & m_2l_1l_2 \cos(x^1-x^2) & l_1 (m_1+m_2) \cos x^1 \\
			m_2l_1l_2 \cos(x^1-x^2) & m_2 l_2^2 & m_2 l_2 \cos x^2 \\
			l_1 (m_1+m_2) \cos x^1 & m_2 l_2 \cos x^2 & m_1+m_2+m_3
		\end{pmatrix}.
	\end{align*} 

The potential energy is 
\begin{align*}
	V=-(m_1+m_2)l_1a\cos x^1-m_2l_2a\cos x^2+  \frac{1}{2}k(x^3)^2,
\end{align*}
where $a$ is the gravitational constant. The equations of the dynamics can be written on $\T Q$, where $Q=\mathbb{S}^1\times \mathbb{S}^1\times \mathbb{R}$, as
\begin{align}
	\label{eq:fdob}
		\begin{split}
		&\dot{x}^1=v^1\qquad
		\dot{v}^1=- \Gamma^1_{11}v^1v^1- \Gamma^1_{22}v^2v^2+ {e}^1+ {g}_1^1  u_1+ {g}_2^1  u_2\\
	&\dot{x}^2=v^2\qquad
	\dot{v}^2=- \Gamma^2_{11}v^1v^1- \Gamma^2_{22}v^2v^2+ {e}^2+ {g}_1^2  u_1+ {g}_2^2  u_2\\
		& \dot{x}^3=v^3\qquad
		 \dot{v}^3=- \Gamma^3_{11}v^1v^1- \Gamma^3_{22}v^2v^2+ {e}^3+ {g}_1^3  u_1+ {g}_2^3  u_2,
	\end{split}
\end{align}

	\begin{align*}
		&M=1/\big(l_1m_1(m_1+m_2)\sin^2x^1+l_1m_2m_3\sin^2(x^1-x^2) +l_1m_1m_3\big)\\
		& \Gamma^1_{11}=\frac{1}{2}Ml_1m_1(m_1+m_2)\sin(2x^1) +\frac{1}{2}Ml_1m_2m_3\sin\left(2(x^1-x^2)\right)\\
		& 	\Gamma^1_{22}=Ml_2m_1m_2 \sin x^1 \cos x^2 + Ml_2 m_2m_3 \sin(x^1-x^2)\\ 
		& 	\Gamma^2_{11}=-\frac{l_1^2(m_1+m_2)m_3}{l_2}M \sin(x^1-x^2)\\
		& 	\Gamma^2_{22}=-\frac{l_1m_2m_3}{2}M\sin\left(2(x^1-x^2) \right) \\
		& 	\Gamma^3_{11}=-l_1^2m_1(m_1+m_2)M \sin x^1\\
		& 	\Gamma^3_{22}=-l_1l_2m_1m_2 M \sin x^1 \cos(x^1-x^2)
			\end{align*}
			\begin{align*}
		&  {g}^1_1=\frac{M}{l_1}\left(m_1+m_3+m_2\sin^2x^2 \right)\\
		&  {g}^2_1= {g}^1_2=-\frac{M}{l_2}\left(m_3\cos(x^1-x^2)+(m_1+m_2)\sin x^1\sin x^2 \right)\\
		&  {g}^3_1=M\left(m_2\sin(x^1-x^2)\sin x^2-m_1\cos x^1 \right)\\ 
		&  {g}^2_2=\frac{l_1(m_1+m_2)}{l_2^2m_2}M\left(m_3+(m_1+m_2)\sin^2x^1 \right) \\
		&  {g}^3_2=-\frac{l_1(m_1+m_2)}{l_2}M \sin x^1 \sin(x^1-x^2)
			\end{align*}
	\begin{align*}
		&  {e}^1=Mm_1(m_1+m_2+m_3)a \sin x^1- Mk x^3 \cos x^1\\&\qquad + Mm_2\sin(x^1-x^2)\left( m_3a \cos x^2+kx^3\sin x^2\right)\\ 
		&  {e}^2=-\frac{l_1(m_1+m_2)}{l_2}M \sin(x^1-x^2)\left( m_3a \cos x^1+k x^3\sin x^1\right)\\
		&  {e}^3=-\frac{l_1m_1(m_1+m_2)a}{2}M\sin(2x^1)+Ml_1m_1kx^3 +Ml_1m_2 k x^3 \sin^2(x^1-x^2).
	\end{align*}

\noindent Again, $e(x)=e^i(x)\fp{}{x^i}$ and $g_r(x)=g_r^i(x)\fp{}{x^i}$, with $e^i(x)$ and $g^i(x)$ defined above, are vector fields on $Q$.

Set the output functions to be $y_1=h_1=x^2$ and $y_2=h_2=x^3$. Immediately we see that $L_{g_r}h_1=g^2_r$ and $L_{g_r}h_2=g^3_r$, for $r=1,2$, therefore the vector relative half-degree is $(\nu_1,\nu_2)=(1,1)$. By applying the following mechanical feedback to system \eqref{eq:fdob} (no change of coordinates needed)
\begin{align*}
	&\ti u_1=- \Gamma^2_{11}v^1v^1- \Gamma^2_{22}v^2v^2+ {e}^2+ {g}_1^2  u_1+ {g}_2^2  u_2\\
	& \ti u_2=- \Gamma^3_{11}v^1v^1- \Gamma^3_{22}v^2v^2+ {e}^3+ {g}_1^3  u_1+ {g}_2^3  u_2,
\end{align*}
we get the system of the form \eqref{eq:md}
\begin{align}
	\label{eq:toras}
	\begin{split}
		&\dot{x}^1=v^1 \qquad 
		\dot{v}^1=- \ti\Gamma^1_{11}v^1v^1- \ti\Gamma^1_{22}v^2v^2+ \ti e^1+  \ti g_1^1\ti u_1+ \ti g_2^1\ti u_2\\
		& \dot{x}^2=v^2\qquad \dot{v}^2=  \ti u_1 \\
		& \dot{x}^3=v^3\qquad\dot{v}^3=\ti u_2
	\end{split}
\end{align}
with the output $(h_1,h_2)=(x^2,x^3)$ and $(x^1,v^1)$-part being unobserved. Denoting $p_1=-\frac{l_2 m_2}{l_1(m_1+m_2)}$ and $p_2=-\frac{m_1+m_2+m_3}{l_1(m_1+m_2)}$, we have
\begin{align*}
	& \ti	\Gamma^1_{11}=-\tan x^1\\
	& \ti	\Gamma^1_{22}=p_1 \sec x^1 \sin x^2\\ 
	&  \ti g^1_1=p_1 \sec x^1 \cos x^2\\
	& \ti g^1_2=p_2 \sec x^1\\
	&  \ti e^1=\frac{k}{l_1(m_1+m_2)}x^3\sec x^1,
\end{align*}
where the vector fields of $\eqref{eq:toras}$ are $\ti e=\ti e^1\fp{}{x^1}$, $\ti g_1=\ti g_1^1 \fp{}{x^1}+\fp{}{x^2}$, and $\ti g_2=\ti g_2^1 \fp{}{x^1}+\fp{}{x^3}$.

Now choose, for system \eqref{eq:toras}, output functions as $ y_1= h_1=\frac{1}{p_1 p_2} \sin x^1-\frac{1}{p_2}\sin x^2 -\frac{x^3}{p_1}$ and $ y_2= h_2=x^2$. A direct calculation shows that $L_{\ti g_1} h_1=L_{\ti g_2} h_1=0$, $L_{\ti e} h_1=p_3 x^3$, for $p_3=\frac{k l_1(m_1+m_2)}{l_2 m_3 (m_1+m_2+m_3)}$, $L_{\ti g_1}L_{\ti e} h_1=0$, $L_{\ti g_2}L_{\ti e} h_1=p_3$. Moreover, we have $L_{\ti g_1}\ti h_2=1$, $L_{\ti g_2}\ti h_2=0$. Now, the vector relative half-degree is $(\ti \nu_1,\ti \nu_2)=(2,1)$ and $\nabla dh_1=0$, therefore, {by Theorem~\ref{thm:MF-out}, the MIOLD-problem is solvable with the unobservable part absent (since $\ti \nu_1+\ti \nu_2 =2+1=n$) and thus, by Corollary \ref{coro:MF-out},} the system is (fully, input-state) MF-linearizable. Indeed, to MF-linearize the system, we need to apply the mechanical diffeomorphism  $\ti x^1=\frac{1}{p_1 p_2} \sin x^1-\frac{1}{p_2}\sin x^2 -\frac{x^3}{p_1}$, $\ti v^1=\frac{1}{p_1 p_2} v^1 \cos x^1-\frac{1}{p_2}v^2\cos x^2 -\frac{v^3}{p_1}$, $\ti x^2=x^2$, $\ti v^2=v^2$, $\ti x^3=x^3$, $\ti v^3=v^3$, {which bring the system into \eqref{eq:md} with $(x_{m+1},v_{m+1})$ absent}. 

\vspace{1cm}

In all above examples we showed that choosing different output functions renders the system partially or fully linearizable. Despite the fact that for mechanical systems the most natural choice for the output functions (available by direct measurements) are the configurations $x^i$ themselves, this could not be the best choice for the control design purposes. It can be more fruitful to use nontrivial output functions $h(x)$ (that MF-linearize the system) and in practical realizations such functions can be computed by the controller (e.g. PC, PLC, etc.).   

\section{Conclusions}

In this paper we have stated and solved the problem of simultaneous input-output feedback linearization and decoupling of mechanical control systems, also known as non-interacting control problem. In contrast to the case of general control systems the problem is nontrivial if we ask for preserving the mechanical structure of the system.
We formulate our conditions in terms of objects defined on the configuration manifold $Q$ only. Therefore a natural question arises, namely to compare the relative degrees $\rho_\ell$ with the newly introduced relative half-degree $\nu_{\ell}$. Furthermore, a related problem is that of the existence of mechanically flat inputs (see preliminary results in \cite{235}), namely the mechanical analogous of the result of \cite{zeitz}. Another interesting topic is to draw a comparison between our result and the problem of the input-output decoupling of Hamiltonian systems, see e.g. \cite{hs1,hs2}. We plan to address all those problems in future works.

\section{Appendix}

The operator of an affine connection $\nabla$ allows to define intrinsically the acceleration as the covariant derivative $\nabla_{\dot{x}(t)}\dot{x}(t)$, see \cite{bullo,bloch,lee}. 
In local coordinates on $Q$, an affine connection is determined by its Christoffel symbols
\begin{align*}
	\nabla_{\fp{}{x^j}}\fp{}{x^k}=\Gamma^i_{jk}\fp{}{x^i},
\end{align*}
and we assume throughout that $\nabla$ is symmetric, i.e. $\Gamma^i_{jk}=\Gamma^i_{kj}$. Given an affine connection $\nabla$, we may differentiate any tensor field on $Q$ along a given vector field, for details see \cite{lee}. We say that a vector field $f$ is \textit{parallel along a curve} $x(t)$ if $ \nabla_{\dot{x}(t)}f=0$. In order to describe a vector field $f=f^k\fp{}{x^k}$ that is parallel along any curve on $Q$, we construct (using linearity of the affine connection) a $(1,1)$-tensor field called \textit{the total covariant derivative} $\nabla f = \nabla_jf^i dx^j\otimes \fp{}{x^i}$, where $\nabla_jf^i=\fp{f^i}{x^j}+\Gamma^i_{jk}f^k$. This object can be viewed as a matrix of covariant derivatives in all directions $\fp{ }{x^j}$. We say that a vector field $f$ is \textit{parallel} on $Q$ if and only if $\nabla f=0$. This concept can be generalized to any tensor field. Below, we give formulae for the total covariant derivatives of functions and differential one-forms, which are used in the paper.

For a $(0,0)$-tensor, i.e. a scalar function $h\in C^{\infty}(Q)$, the total covariant derivative $\nabla h$ is the exact differential $dh\in \Lambda(Q)$, where $\Lambda(Q)$ is the set of smooth differential one-forms, that is the $(0,1)$-tensor
\begin{align*}
	\nabla h =dh=\fp{h}{x^j}dx^j.
\end{align*} 

For a $(0,1)$-tensor field, i.e. a differential one-form $\omega={\omega_j\ dx^j}\in \Lambda(Q)$, the total covariant derivative $\nabla\omega=\nabla_k\omega_j \ dx^j\otimes dx^k$ is the $(0,2)$-tensor~field, with
\begin{align*}
	\nabla_k\omega_j=\fp{\omega_j}{x^k}-\Gamma^i_{jk}\omega_i.
\end{align*} 
Finally, by a combination of the above formulae, we have
\begin{align}
	\label{eq:cdwa}
 \nabla dh=(\fp{^2 h}{x^j\partial x^k}-\Gamma^i_{jk}\fp{h}{x^i})\ dx^j\otimes dx^k.
\end{align}

\subsection*{Acknowledgements}
We are grateful to anonymous Reviewers for their constructive
criticism and suggestions that helped to improve the final
presentation of the paper.\\

This work was partially supported by the Statutory Grant No. 0211/SBAD/0123, 0211/SBAD/0324.

\end{document}